\documentclass[12pt]{amsart}
\usepackage{enumerate,amssymb,version,geometry,aliascnt,bookmark}
\usepackage[all]{xy}
\usepackage{hyperref}

\hypersetup{
pdfauthor={Olivier Haution},
pdftitle={Integrality of the Chern character in small codimension},
pdfkeywords={Steenrod operations, Chow groups, Chern character, Grothendieck-Riemann-Roch theorem, positive characteristic},
hidelinks,
pdftex,
linktocpage=false
}

\DeclareMathOperator{\id}{id}
\DeclareMathOperator{\Spec}{Spec}
\DeclareMathOperator{\rank}{rank}
\DeclareMathOperator{\codim}{codim}
\DeclareMathOperator{\trdeg}{tr.deg.}
\DeclareMathOperator{\CH}{CH}
\DeclareMathOperator{\Ch}{Ch}
\DeclareMathOperator{\CHc}{\mathsf{CH}}
\DeclareMathOperator{\ch}{ch}
\DeclareMathOperator{\chc}{\mathsf{ch}}
\DeclareMathOperator{\gr}{gr}
\newcommand{\colim}[1]{\mycolim_{#1}}
\DeclareMathOperator*{\mycolim}{colim}

\DeclareMathOperator{\im}{im}
\DeclareMathOperator{\Sq}{Sq}
\DeclareMathOperator{\Td}{Td}

\newcommand{\Tan}{T}
\newcommand{\Var}{\mathsf{Var}}
\newcommand{\Varp}{\mathsf{Proj}}
\newcommand{\Ab}{\mathsf{Ab}}
\newcommand{\Zz}{\mathbb{Z}}
\newcommand{\Qq}{\mathbb{Q}}
\newcommand{\CHQ}[1]{\CH({#1})_{\Qq}}
\newcommand{\CHQn}[2]{\CH_{#2}({#1})_{\Qq}}
\newcommand{\CHt}{\widetilde{\CH}}

\newcommand{\CHZ}[1]{\CH({#1})_{\Zz\subset \Qq}}
\newcommand{\CHZpn}[2]{\CH_{#1}({#2})_{\Zz_{(p)}\subset \Qq}}
\newcommand{\CHZn}[2]{\CH_{#1}({#2})_{\Zz\subset \Qq}}
\newcommand{\nx}[1]{n_{#1}(p)}
\newcommand{\xx}{\{x\}}

\newcommand{\CHcQ}[1]{\CHc({#1})_{\Qq}}

\newcommand{\G}{K_0'}
\newcommand{\K}{K_0}

\newcommand{\Oc}{\mathcal{O}}

\newcommand{\para}[1]{\paragraph { \emph{#1}.}}
\newcommand{\Te}{T'}

\newtheorem{theorem}{Theorem}[section]

\newaliascnt{proposition}{theorem}
\newtheorem{proposition}[proposition]{Proposition}
\aliascntresetthe{proposition}

\newaliascnt{lemma}{theorem}
\newtheorem{lemma}[lemma]{Lemma}
\aliascntresetthe{lemma}

\newaliascnt{statement}{theorem}
\newtheorem{statement}[statement]{Statement}
\aliascntresetthe{statement}

\newaliascnt{corollary}{theorem}
\newtheorem{corollary}[corollary]{Corollary}
\aliascntresetthe{corollary}

\newtheorem*{conjecture}{Conjecture}

\theoremstyle{definition}

\newaliascnt{remark}{theorem}
\newtheorem{remark}[remark]{Remark}
\aliascntresetthe{remark}

\newaliascnt{example}{theorem}

\aliascntresetthe{example}

\newaliascnt{definition}{theorem}

\aliascntresetthe{definition}

\begin{document}
\begin{abstract}
We prove an integrality property of the Chern character with values in Chow groups. As a consequence we obtain a construction of the $p-1$ first homological Steenrod operations on Chow groups modulo $p$ and $p$-primary torsion, over an arbitrary field. We provide applications to the study of correspondences between algebraic varieties.
\end{abstract}

\author{Olivier Haution}
\title[Integrality of the Chern character]{Integrality of the Chern character in small codimension}
\email{olivier.haution at gmail.com}
\address{School of Mathematical Sciences\\
University of Nottingham\\
University Park\\
Nottingham\\
NG7 2RD\\ 
United Kingdom}
\date{\today}


\maketitle

\section*{Introduction}
The Grothendieck Riemann-Roch theorem of \cite{BFM} can be expressed as the existence of the \emph{homological Chern character}, a collection of group morphisms
\[
\ch_{i} \colon \G(X) \to \CHQn{X}{i},
\]
which are compatible with projective push-forwards. Here $X$ is a (possibly singular) algebraic variety, $\G(X)$ the Grothendieck group of coherent sheaves, and $\CHQn{X}{i}$ the Chow group of $X$ of algebraic cycles of dimension $i$, with $\Qq$-coefficients. Compatibility with pull-backs is expressed by a Riemann-Roch type formula involving the Todd class. 

Given a prime number $p$, one can formulate the following statement describing a $p$-integrality property of the homological Chern character.
\begin{conjecture}
\label{conjp}
Let $X$ be a variety, and $n$ a positive integer. Then
\[
p^{[n/(p-1)]}\cdot \ch_{\dim X-n}[\Oc_X] \in \im(\CH_{\dim X-n}(X)_{\Zz_{(p)}} \to \CH_{\dim X-n}(X)_{\Qq}).
\]
\end{conjecture}
The topological analog of this conjecture was stated and proved by Adams in \cite[Theorem~1]{Ada-On-61}. A more elementary proof of this topological result was later given by Atiyah (\cite[Theorem~7.1]{Ati-Po-66}), under the quite restrictive assumption that $X$ is a CW-complex whose homology has no torsion. Although our approach is more akin to Atiyah's than Adams's (where in particular Steenrod operations were used), no condition of such nature will appear.

In \cite[Theorem~4.4]{firstsq}, we proved the conjecture for $p=2$ and $n=1$. In the present article, we prove it for $n < p(p-1)$, and any $p$ (\autoref{sect:unconditional}). In particular we recover the previously obtained result, using an argument of different nature.
 We discuss in \autoref{sect:base} how this result can be extended to the more general setting of a scheme of finite type over a regular base.

Let us mention that the conjecture is an immediate consequence of the formula expressing compatibility of the homological Chern character with pull-backs when $X$ is smooth, or more generally a local complete intersection variety. This remark can be used to prove the conjecture in the presence of some form of resolution of singularities. Using a result of Gabber on existence of regular alterations, we prove the conjecture for all $n$, when the characteristic of the base field is not $p$.\\

The $n$-th Steenrod operation on Chow groups modulo $p$ is an endomorphism of the Chow group with $\Zz/p$-coefficients which lowers dimension of cycles by $n(p-1)$. The Steenrod operations constitute an efficient tool for studying rationality properties of projective homogeneous varieties.

The classical constructions of the operations modulo $p$ (\cite{Boi-A-08}, \cite{Bro-St-03}, \cite{EKM}, \cite{Lev-St-05}, \cite{Vo-03}) do not work over a field of characteristic $p$. In \cite{firstsq}, we constructed a weak form of the operation for $p=2$ and $n=1$, over any field. We later obtained the full version of this operation in \cite{duality}. In the present paper, we construct a weak form of the $n$-th Steenrod operation modulo $p$, for $p$ an arbitrary prime and $n=1,\cdots,p-1$, over any field. 

More precisely, assuming that the conjecture above holds true for $n\leq m(p-1)$, and given any (possibly singular) variety $X$ over an arbitrary field, we construct for each integer $i$ such that $0<i\leq m$, a morphism of graded abelian groups
\[
T_i \colon \CH_\bullet(X) \otimes \Zz/p \to \CHt_{\bullet-i(p-1)}(X) \otimes \Zz/p.
\]
Here $\CHt_d$ is the quotient of the Chow group $\CH_d$ by its torsion subgroup. These morphisms are compatible with proper push-forwards, external products, and extension of the base field. We also provide a Riemann-Roch type formula expressing compatibility with pull-backs along local complete intersection morphisms.

Thus we construct operations $T_1,\cdots,T_{p-1}$ unconditionally, and operations $T_i$ for all $i>0$ when the characteristic of the base field is not $p$. The operation $T_i$ corresponds --- up to a sign --- to the $i$-th Steenrod operation modulo $p$ when $i\leq p$. More generally the total operation $\sum_i T_i$ corresponds to the inverse of the total Steenrod operation $\sum_j S_j$. We explain how to recover the operation induced by $S_j$ on reduced Chow groups (the \emph{reduced Steenrod operation}) from the operations $T_i$.\\

In the last part of this paper, we prove the degree formula corresponding to the operations $T_i$. Since we only proved the conjecture in general for small values of $n$, we need to work with varieties of dimension $\leq p(p-1)$ (note that the dimension $p(p-1)$ is allowed). We first give a bound on the $p$-adic valuation of the index of a projective variety in terms of its dimension and of the Euler characteristic of its structure sheaf. Then we deduce from the degree formula that the Severi-Brauer variety of a central division algebra of prime degree $p$ is strongly $p$-incompressible. This is a new statement when the base field has characteristic $p$. Finally we show that the Euler characteristic of the structure sheaf of a smooth projective variety of dimension $p-1$, possessing a special correspondence but no zero-cycle of degree prime to $p$, is prime to $p$. This extends a statement of Rost to arbitrary base fields.

\setcounter{tocdepth}{1}
\tableofcontents

\section{Notations}
\para{Varieties} We fix a base field $k$. A \emph{variety} is a finite type, separated, quasi-projective scheme over $k$. A morphism of varieties is a morphism of schemes over $k$. We will denote by $\Var/k$ this category. We also denote by $\Varp/k$ the category with the same objects, but with morphisms only those which are projective. The function field of an integral variety $X$ shall be denoted by $k(X)$.\\

\para{Integral elements} Given a rational number $\alpha$, we write $[\alpha]$ for the greatest integer $\leq \alpha$. 

Let $A$ be an abelian group. For any commutative ring $R$ we write $A_R$ for $A \otimes_{\Zz} R$, and, if $S$ is another commutative ring containing $R$, we write
\[
A_{R \subset S} = \im (A_R \to A_S).
\]
Note that if $A$ is a ring, then $A_{R \subset S}$ is a subring of $A_S$.

Let $p$ be a prime number. When $R=\Zz_{(p)}$ (the subring of $\Qq$ consisting of fractions with denominator prime to $p$), $S=\Qq$, and $x \in A_{\Qq}-\{0\}$, we denote by $v_p(x)$ the greatest integer $r$ such that $p^{-r}\cdot x \in A_{\Zz_{(p)} \subset \Qq}$. We set $v_p(0)=\infty$. 

When $x,y \in A_{\Qq}$, we have
\[
v_p(x+y) \geq \min(v_p(x),v_p(y)).
\]

If $A$ is commutative ring, $M$ an $A$-module, $a \in A_{\Qq}$, and $m \in M_{\Qq}$, we have
\begin{equation}
\label{eq:vpmod}
v_p(a \cdot m)\geq v_p(a) + v_p(m).
\end{equation}\\

\para{Grothendieck groups of schemes} 
We denote by $\K(-)$ the presheaf of rings on $\Var/k$ which associates to a variety $X$ the Grothendieck ring of locally-free coherent $\Oc_X$-modules. We denote by $\G(-)$ the functor $\Varp/k \to \Ab$ which associates to a variety $X$ the Grothendieck group of coherent $\Oc_X$-modules. Here, and in the rest of the paper, the notation $\Ab$ stands for the category of abelian groups.\\

\para{Local complete intersection morphisms} (See \cite[Appendix~B.7.6]{Ful-In-98}) These are the morphisms of varieties $f\colon Y \to X$ which can factored as $p\circ i$, with $i$ a regular closed embedding, and $p$ a smooth morphism. Such a morphism admits a \emph{virtual tangent bundle} $\Tan_f=i^*[\Tan_p]-[N_i] \in \K(Y)$, where $\Tan_p$ is the tangent bundle of $p$, and $N_i$ the normal bundle of $i$. This element does not depend on the choice of the factorization. There are pull-backs $f^* \colon \CH(X) \to \CH(Y)$ and $f^*\colon \G(X) \to \G(Y)$. We call a variety $X$ a \emph{local complete intersection variety} if its structural morphism $x\colon X \to \Spec(k)$ is a local complete intersection morphism. We denote by $\Tan_X \in \K(X)$ the virtual tangent bundle $\Tan_x$ of $x$.\\

\para{Topological filtration}
Let $X$ be a variety, and $d$ an integer. We denote by $\G(X)_{(d)}$ the subgroup of $\G(X)$ generated by those $\Oc_X$-modules whose support has dimension $\leq d$. It coincides with the subgroup generated by elements of type $i_*[\Oc_Z]$, where $i\colon Z \hookrightarrow X$ is a closed embedding, and $\dim Z \leq d$. 

When $f \colon Y \to X$ is a projective morphism, we have 
\[
f_*\G(Y)_{(d)} \subset \G(X)_{(d)}.
\]

When  $f \colon Y \to X$ is a local complete intersection morphism of relative dimension $n$, we have (\cite[Theorem~83 and Lemma~84]{Gil-K-05})
\[
f^*\G(X)_{(d)} \subset \G(Y)_{(d+n)}.
\]

The associated graded group is denoted by
\[
\gr_d\G(X)=\G(X)_{(d)}/\G(X)_{(d-1)}.
\]
There is a natural transformation of functors $\Varp/k \to \Ab$
\begin{equation}
\label{eq:chgrk}
\varphi_{-}\colon \CH(-) \to \gr \G(-)
\end{equation}
which respects the gradings (\cite[Example~15.5]{Ful-In-98}). It is moreover compatible with pull-backs along local complete intersection morphisms, external products, and extension of the base field (see \cite[Section~3]{reduced}).\\

\para{Chow cohomology ring}
Let $X$ be a possibly singular variety, and consider the operational Chow ring 
$\CH(X \xrightarrow{\id} X)$ (see \cite[Definition~17.3]{Ful-In-98}). When resolution of singularities is available, this ring is commutative (\cite[Example~17.4.4]{Ful-In-98}). We denote by $\CHc(X)$ the center of $\CH(X \xrightarrow{\id} X)$. The group $\CH(X)$ has a natural structure of a $\CHc(X)$-module. We obtain a presheaf of commutative rings $\CHc$ on $\Var/k$, in which Chern classes of vector bundles take their values (\cite[Proposition~17.3.2]{Ful-In-98}).  

When $E$ is a vector bundle over a variety $X$, we write
\[
\Td(E) \in \CHc(X)_{\Qq}
\]
for the total Todd class (\cite[Example~3.2.4]{Ful-In-98}). The individual components are
\[
\Td^i(E) \in \CHc^i(X)_{\Qq}.
\]\\

\para{Chern character} The homological Chern character is a natural transformation of functors $\Varp \to \Ab$  
\[
\ch \colon \G(-) \to \CHQ{-}.
\]
This is the map $\tau$ of \cite[Theorem 18.3]{Ful-In-98}. Its component of dimension $i$ will be
\[
\ch_i \colon \G(-) \to \CH_i(-)_{\Qq}.
\]\\

\para{Adams operations} (See \cite{Sou-Op-85}) Let $l \in \Zz - \{0\}$. The $l$-th homological Adams operation is a natural transformation of functors $\Varp \to \Ab$
\[
\psi_l \colon \G(-) \to \G(-) \otimes \Zz[1/l] .
\]

\part{Integrality of the Chern character}
We fix an prime number $p$. We say that \emph{$p$-integrality holds in codimension $\leq c$} if the following statement is true.
\begin{statement}
\label{st:statement}
Let $X$ be a variety and $d\geq 0$ an integer. Then for all integers $n$ such that $0\leq n \leq c$, and all elements $x \in \G(X)_{(d)}$, we have 
\[
v_p(\ch_{d -n}(x)) \geq -\Big[\frac{n}{p-1}\Big].
\]
\end{statement}
We prove in \autoref{sect:unconditional} that $p$-integrality holds in codimension $\leq p(p-1)-1$. Moreover, we prove in  \autoref{sect:gabber} that $p$-integrality holds in codimension $\leq m$ for all $m$ when the characteristic of the base field is not $p$. We expect that this assumption is superfluous.

\section{Adams operations and Chern character}
\label{sect:adamschern}
We fix an element $l \in \Zz- \{0\}$. The sole purpose of this section is to prove, in \autoref{prop:chpsi}, the relation for $n\geq 0$,
\begin{equation}
\label{eq:chpsi}
\ch_n \circ \psi_l=l^{-n} \cdot \ch_n.
\end{equation}
This is certainly very classical, at least so is the cohomological counterpart of this statement. Due to the lack of suitable reference, we give some details in this section. We shall need to introduce some notations, but except in \autoref{sect:base}, the only thing used in the sequel will be the relation \eqref{eq:chpsi}.\\

\para{Cohomological variants} The cohomological (i.e. usual) Chern character (\cite[Example~3.2.3]{Ful-In-98}) is a morphism of presheaves of rings on $\Var/k$
\[
\chc \colon \K(-) \to \CHcQ{-}.
\]

The cohomological (i.e. usual) $l$-th Adams operation  is a morphism of presheaves of rings on $\Var/k$
\[
\psi^l \colon \K(-) \to \K(-),
\]
sending the class of a line bundle to its $l$-th tensor power.

When $X$ is connected, the ring $\K(X)$ has an augmentation
\[
\rank \colon \K(X) \to \Zz
\]
sending the class of a vector bundle to its rank.\\

\para{Bott's class} Consider the Laurent polynomial in the variable $x$
\[
t^l(x)=\frac{1-x^l}{1-x}.
\]

There is a unique morphism of presheaves of abelian groups on $\Var/k$
\[
\theta^l \colon \K(-) \to \big(\Zz[1/l] \otimes \K(-)\big)^{\times}
\]
such that for all line bundles $L$,
\[
\theta^l(L)= t^l[L^\vee].
\]
Note that this does not follow \cite{Sou-Op-85}, where $\theta^{l}(L)=t^l[L]$.\\

\para{$K$-theory with supports} When $X$ is a closed subvariety of a variety $W$, the notation $\K^X(W)$ stands for the Grothendieck group of complexes of locally free $\Oc_W$-modules, which are acyclic off $X$, modulo the classes of acyclic complexes. It has a natural structure of a $\K(W)$-module.  We denote by $\CH(X\to W)$ the bivariant Chow group (see \cite[Definition~17.1]{Ful-In-98}). There are maps
\[
-\cap [\Oc_W] \colon \K^X(W) \to \G(X),\quad -\cap [W] \colon \CH(X\to W) \to \CH(X)
\]
which are isomorphisms when $W$ is smooth (\cite[Lemma~1.9]{GS-87} and \cite[Proposition~17.4.2]{Ful-In-98}). We write
\[
\chc_X^W \colon \K^X(W) \to \CHQ{X\to W}.
\]
for the local Chern character (\cite[\S 18.1]{Ful-In-98}), and
\[
\psi^l_X \colon \K^X(W) \to \K^X(W)
\]
for the Adams operation with supports (\cite[Section~4]{GS-87}). 

Assume that $W$ is smooth, and that $x \in \K^X(W)$. Then the homological Chern character is given by the formula (\cite[Theorem~18.3, (3)]{Ful-In-98}),
\begin{equation}
\label{eq:homch}
\ch(x\cap [\Oc_W])=(\Td(\Tan_W) \cdot \chc_X^W(x)) \cap [W],
\end{equation}
while homological Adams operations are defined by the formula (\cite[Th\'eor\`eme~7]{Sou-Op-85})
\begin{equation}
\label{eq:homad}
\psi_l(x \cap [\Oc_W])=(\theta^l(-\Tan_W) \cdot \psi^l_X(x)) \cap [\Oc_W].
\end{equation}

Composing the maps $\chc_X^W$ with projections  $\CHQ{X\to W} \to \CH^j(X \to W)_{\Qq}$, we obtain individual components $(\chc_X^W)^j$.\\

We now prove a few lemmas that will be used to obtain \eqref{eq:chpsi}.
\begin{lemma}
\label{lemm:tdpsi}
Let $X$ be a connected variety. We have in $\CHcQ{X}$, for all $u \in \K(X)$, 
\[
\Td(-u) \cdot \chc \circ \theta^l(u)=l^{\rank(u)} \cdot \Td(-\psi^l(u)).
\]
\end{lemma}
\begin{proof}
Assume that $L$ is a line bundle. Then  
\[
\chc \circ \theta^l(L)=\chc(t^l[L^\vee])=t^l(\chc[L^\vee])=t^l(x^{-1})=\frac{1-x^{-l}}{1-x^{-1}},
\]
where $x=e^{c_1(L)}$. Then  
\begin{align*} 
 \Td(-L) \cdot \chc \circ \theta^l(L) &=  \frac{1-x^{-1}}{c_1(L)} \cdot \frac{1-x^{-l}}{1-x^{-1}}\\ 
  &=\frac{1-x^{-l}}{c_1(L)}\\
  &=l\cdot \frac{1-e^{-c_1(L^{\otimes l})}}{c_1(L^{\otimes l})}\\
  &=l\cdot \Td(-L^{\otimes l}).
   \end{align*}
Thus the claimed formula holds when $u$ is the class of a line bundle. Both sides of this formula define morphisms of presheaves of abelian groups on $\Var/k$
\[
\K(-) \to (\CHc(-)_{\Qq})^{\times}.
\]
Therefore we can conclude with the splitting principle.
\end{proof}

\begin{lemma}
\label{lemm:tdpsi2}
Let $X$ be a variety. We have in $\CHcQ{X}$, for all $u \in \K(X)$, 
\[
\Td(\psi^l(-u))=\sum_{i\geq 0} l^i \cdot \Td^i(-u).
\]
\end{lemma}
\begin{proof}
For a line bundle $L$, we have 
\[
\Td(-L^{\otimes l})=\sum_{i \geq 0} \frac{(-c_1(L^{\otimes l}))^i}{(i+1)!}=\sum_{i \geq 0} l^i \cdot \frac{(-c_1(L))^i}{(i+1)!}=\sum_{i \geq 0} l^i \cdot \Td^i(-L),
\]
proving the formula for $u$ the class of a line bundle. The lemma follows as above from the splitting principle.
\end{proof}

\begin{lemma}[{\cite[Theorem~3.1]{Kurano-Roberts}}]
\label{lemm:tdpsi3}
Let $X$ be a closed subvariety of a smooth variety $W$. We have in $\CHQ{X \to W}$, for $x \in \K^X(W)$,
\[
\chc_X^W(\psi^l_X(x))=\sum_{j\geq 0} l^j \cdot (\chc_X^W)^j(x).
\]
\end{lemma}

\begin{proposition}
\label{prop:chpsi}
Let $n\geq 0$ be an integer. We have $\ch_n \circ \psi_l=l^{-n} \cdot \ch_n$.
\end{proposition}
\begin{proof}
It will be enough to prove that the two sides of the formula have the same effect on any element of $\G(X)$, for $X$ a connected variety. Choose a closed embedding $X \hookrightarrow W$ into a smooth connected variety $W$ of dimension $d$. Then, we have in $\CHQ{X}$, for all $x \in \K^X(W)$,
\begin{align*} 
 \ch \circ \psi_l(x\cap [\Oc_W]) &= \Td(\Tan_W|_X)\cdot \big( \chc_X^W(\theta^l(-\Tan_W) \cdot \psi^l_X(x)) \cap [W]\big) \\ 
  &= \big( \Td(\Tan_W|_X) \cdot \chc \circ \theta^l(-\Tan_W|_X) \big) \cdot \big(\chc_X^W(\psi^l_X(x))\cap [W]\big)\\
  &= l^{-d} \cdot \Td(\psi^l(\Tan_W|_X))\cdot \big(\chc_X^W(\psi^l_X(x))\cap[W]\big)\\
  &=\sum_{i=0}^d l^{i-d} \cdot \Td^i(\Tan_W|_X) \cdot \sum_{j\geq 0} l^j\cdot \big((\chc_X^W)^j(x)\cap [W]\big)\\
  &=\sum_{n=0}^d l^{-n} \cdot \sum_{i+j=d-n} \Td^i(\Tan_W|_X) \cdot \big((\chc_X^W)^j(x)\cap [W]\big)\\
  &=\sum_{n=0}^d l^{-n} \cdot \ch_n(x\cap[\Oc_W]).
   \end{align*}
The first equality is \eqref{eq:homch} and \eqref{eq:homad} . The second follows from \cite[Proposition~18.1, (c)]{Ful-In-98}. For the third one, we use \autoref{lemm:tdpsi}. The fourth one uses \autoref{lemm:tdpsi2} and \autoref{lemm:tdpsi3}. The last one is again \eqref{eq:homch}.
\end{proof}

\section{Integrality in arbitrary characteristic}
\label{sect:unconditional}

\begin{proposition}
\label{prop:graded}
Let $l \in \Zz - \{0\}$, $X$ a variety, and $x \in \G(X)_{(d)}$. Then
\[
\psi_l(x)-l^{-d} \cdot x \in \G(X)_{(d-1)} \otimes \Zz[1/l].
\]
\end{proposition}
\begin{proof}
When $l=p$ is a prime number, this is \cite[Lemma~6.1]{reduced}. Actually we obtained this result using the formula 
\[
\theta^p(-d)=p^{-d}, \quad (d \in \Zz), 
\]
which is true when $p$ is an arbitrary element of $\Zz - \{0\}$. Therefore we can obtain the proposition using the proof of \cite[Lemma~6.1]{reduced}.

Alternatively, one can assume that $X$ is integral and $d=\dim X$. Then we have an exact sequence
\[
0 \to \G(X)_{(d-1)} \otimes \Zz[1/l] \to \G(X)_{(d)}\otimes \Zz[1/l] \to \G(\Spec(k(X)))\otimes \Zz[1/l] \to 0,
\]
and the statement is a consequence of \cite[Lemma~4.2]{duality}. 
\end{proof}

\begin{theorem}
\label{th:mainvp}
$p$-integrality holds in codimension $\leq p(p-1)-1$.
\end{theorem}
\begin{proof}
We must prove \autoref{st:statement} with $c=p(p-1)-1$. For all $d$, the case $n=0$ immediately follows from \cite[Theorem~18.3, (5)]{Ful-In-98}. In particular the theorem is proved when $d=0$. We can therefore proceed by induction on $d$, and assume that $n >0$. 

For all $i\geq 1$, define an integer 
\[
a_p(i)= \left\{ \begin{array}{rl}
		  1 &\mbox{ if $p-1|i$} \\
		  0 &\mbox{ otherwise}
       		\end{array}
	\right.
\]
so that we have
\begin{equation}
\label{eq:sumap}
\Big[\frac{n}{p-1}\Big]=\Big[\frac{n-1}{p-1}\Big] + a_p(n).
\end{equation}

Since $0 <n<p(p-1)$, we can use \autoref{lemm:pl} below to choose a integer $l$, prime to $p$, such that
\begin{equation}
\label{eq:vpln}
v_p(l^n-1)=a_p(n).
\end{equation}

Write, in $\G(X) \otimes \Zz[1/l]$, using \autoref{prop:graded}, 
\[
\psi_l(x)=l^{-d}\cdot x + l^{-d-m} \cdot \alpha 
\]
with $\alpha \in \G(X)_{(d-1)}$, and $m\geq 0$. Then applying $l^{d+m} \cdot \ch_{d-n}$ to the relation above, and using \autoref{prop:chpsi}, we obtain in $\CH_{d-n}(X)_{\Qq}$
\[
l^m\cdot (l^n-1) \cdot \ch_{d-n}(x)= \ch_{d-n}(\alpha).
\]
Because of \eqref{eq:vpln}, we get
\[
v_p(\ch_{d-n}(x))=v_p(\ch_{d-n}(\alpha)) + a_p(n). 
\]

We conclude using \eqref{eq:sumap} and the induction hypothesis for $\alpha \in \G(X)_{(d-1)}$.
\end{proof}

\begin{lemma}
\label{lemm:pl}
Let $p$ be a prime number. Then there exists an integer $l$, prime to $p$, with the following properties:
\begin{enumerate}[(a)]
\item If $n$ is not divisible by $p-1$, then $l^n-1$ is not divisible by $p$.

\item If $n$ is not divisible by $p(p-1)$, then $l^n-1$ is not divisible by $p^2$.
\end{enumerate}
\end{lemma}
\begin{proof}
Take for $l$ an integer whose class modulo $p^2$ generates the cyclic group $(\Zz/p^2)^{\times}$. Then one checks easily that the class of $l$ modulo $p$ generates the cyclic group $(\Zz/p)^{\times}$.
\end{proof}

\begin{remark}
In contrast with \cite{firstsq} where the normalization procedure was used, we do not use geometric arguments here. We only use basic relations between Adams operations, Chern character, and filtration on $K$-theory. Using \cite[(1.5)]{Adams-Atiyah} in place of \autoref{prop:graded}, and \cite[Theorem~5.1, (vi)]{Adams-Vector} in place of \autoref{prop:chpsi}, one could proceed as above to replace the assumption ``$X$ is without torsion'' by ``$n<p(p-1)$'' in \cite[Theorem~7.1]{Ati-Po-66}, thus giving an elementary proof of the $p$-primary part of \cite[Theorem~1]{Ada-On-61}, for small $n$. 

Note also that the approach of \cite{firstsq} is limited by nature to schemes whose normalization morphism is of finite type, while \autoref{th:mainvp} is valid for more general schemes (see \autoref{sect:base}).
\end{remark}

\section{Consequences of Gabber's Theorem}
\label{sect:gabber}
In this section, we use the theorem of Gabber asserting that, over a field of characteristic different from $p$, any integral variety $X$ admits a regular alteration of degree prime to $p$. This means that there exists a projective morphism $Z \to X$, generically finite of degree prime to $p$, with $Z$ regular. Some details can be found in \cite{Gabber} and \cite{Gabber-Illusie}.\\

If $d\geq 0$ is an integer, the $d$-th \emph{Todd number} is
\begin{equation}
\label{Todd}
\tau_d=\prod_{p \text{ prime}} p^{[d/(p-1)]}.
\end{equation}

\begin{proposition}
\label{prop:smooth2}
Let $X$ be a local complete intersection variety. Then, for all $n$
\[
\tau_n\cdot \ch_{\dim X -n}[\Oc_X] \in \CH(X)_{\Zz \subset \Qq}.
\]
\end{proposition}
\begin{proof}
Let $x \colon X \to \Spec(k)$ be the structural morphism of $X$. Then by \cite[Theorem~18.3, (4) and (5)]{Ful-In-98},
\[
\ch[\Oc_X]=\ch\circ x^*[\Oc_{\Spec(k)}]=\Td(\Tan_X)\cdot x^*\circ \ch[\Oc_{\Spec(k)}]=\Td(\Tan_X)\cdot[X].
\]
Then the proposition follows from \autoref{lemm:toddvp} below.
\end{proof}

\begin{theorem}
\label{th:mainvl}
Let $p$ be a prime different from the characteristic of the base field. Then $p$-integrality holds in codimension $\leq c$, for any $c$.
\end{theorem}
\begin{proof}
We must prove \autoref{st:statement} with $c$ arbitrary large. We proceed by induction on $d$, the case $d=0$ being obvious. It is enough to consider the case $X$ integral of dimension $d$, and $x=[\Oc_X]$. Using Gabber's theorem, we choose $f \colon Z \to X$ a projective morphism, with $Z$ a regular connected variety of dimension $d$, and such that $[k(Z):k(X)]=\lambda$ is prime to $p$. Then there is an element $\delta \in \G(X)_{(d-1)}$ such that
\[
f_*[\Oc_Z]=\lambda \cdot x + \delta.
\]
Thus
\begin{align*}
v_p(\ch_{d -n}(x))&\geq \min\big(v_p(f_*\ch_{d -n}[\Oc_Z]) , v_p(\ch_{d -n}(\delta))\big)\\
&\geq \min\big(v_p(\ch_{d -n}[\Oc_Z]) , v_p(\ch_{(d-1) -(n-1)}(\delta))\big).
\end{align*}
We conclude using \autoref{prop:smooth2} for $Z$, and the induction hypothesis on $\delta$.
\end{proof}

The combination of \autoref{th:mainvp} and \autoref{th:mainvl} implies:
\begin{corollary}
\label{th:main}
Let $X$ be a variety over a field of characteristic $p$, and $d$ be an integer. Then for all integers $n$ such that $0\leq n < p(p-1)$, and all elements $x \in \G(X)_{(d)}$, we have 
\[
\tau_n \cdot \ch_{d -n}(x) \in \CH(X)_{\Zz \subset \Qq}.
\]
\end{corollary}

\part{Operations on Chow groups}
\section{Construction}
\label{sect:construction}
We proceed as in \cite{firstsq} to derive operations from the integrality property of the Chern character stated above. This follows \cite[Section~7]{Ati-Po-66}.\\

Let $A$ be an abelian group. Then the natural map
\begin{equation}
\label{eq:cht}
A_{\Zz \subset \Qq} \otimes \Zz/p \to A_{\Zz_{(p)} \subset \Qq} \otimes \Zz/p, 
\end{equation}
is an isomorphism. This group can also be described as $A$ modulo $p$ and $p$-primary torsion. By this we mean the quotient of $A$ by $p\cdot A+T$, where $T$ is the subgroup of elements $a \in A$ such that $p^n\cdot a=0$ for some integer $n>0$.\\ 

We assume that $p$-integrality holds in codimension $\leq m(p-1)$, for some integer $m$. Let $q>0$ be an integer. For every integer $i$ such that $0<i \leq m$, we see using \autoref{st:statement}, for $d=q$ and $n=i(p-1)$, that the association $x \mapsto p^i \cdot \ch_{q-i(p-1)}(x)$ induces a natural transformation of functors $\Varp \to \Ab$
\[
\G(-)_{(q)} \to \CHZpn{q-i(p-1)}{-}.
\]
Let $X$ be a variety, and $x \in \G(X)_{(q-1)} \subset  \G(X)_{(q)}$. Then 
\[
p^i \cdot \ch_{q-i(p-1)}(x)=p \cdot (p^{i-1} \cdot \ch_{(q-1)-(i(p-1)-1)}(x)).
\]
Using \autoref{st:statement} for  $d=q-1$ and $n=i(p-1)-1$, we see that this element belongs to $p \cdot \CHZpn{q-i(p-1)}{X}$. This gives a natural transformation of functors $\Varp \to \Ab$
\[
\gr_q \G(-) \otimes \Zz/p \to \CHZpn{q-i(p-1)}{-} \otimes \Zz/p.
\]

Composing on the left with the map \eqref{eq:chgrk}, and on the right with the inverse of the map \eqref{eq:cht} with $A=\CH_{q-i(p-1)}(-)$, we obtain a natural transformation of functors $\Varp \to \Ab$
\[
T_i\colon \CH_q(-) \otimes \Zz/p \to \CHZn{q-i(p-1)}{-} \otimes \Zz/p.
\]
It will be convenient to agree that 
\[
T_0\colon \CH_q(-) \otimes \Zz/p \to \CHZn{q}{-} \otimes \Zz/p
\]
is the quotient map modulo the image modulo $p$ of $p$-primary torsion.

\section{The inverse Todd genus}
We consider the power series in the variable $x$
\begin{equation}
\label{eq:r}
R^{(p)}(x)=\sum_{i=0}^\infty (-1)^i \cdot x^{p^i-1}.
\end{equation}

There is a unique morphism of presheaves of abelian groups on $\Var/k$
\[
r^{(p)} \colon \K(-) \to (\CHc(-)\otimes \Zz/p)^{\times}
\]
such that for $L$ a line bundle, 
\[
r^{(p)}(L)=R^{(p)}(c_1(L)).
\]
Individual components are denoted
\[
r^{(p)}_i \colon \K(-) \to \CHc^{i(p-1)}(-)\otimes \Zz/p.
\]
Using the terminology of \cite{Pan-Ri-02}, this class will be the \emph{inverse Todd genus} of the operation $T=\sum_i T_i$ constructed in the previous section. In this section we describe a relation between the characteristic class $r^{(p)}$ and the Todd class.

\begin{remark}
\label{rem:inv}
The inverse Todd genus of the usual Steenrod operation $\sum_i S_i$ is defined by the power series $W^{(p)}(x)=1+x^{p-1}$. One easily checks that $x\cdot R^{(p)}(x)$ is the left inverse of $x\cdot W^{(p)}(x)$ for the law given by composition of power series with coefficients in $\Zz/p$. 
\end{remark}

We begin by stating an easy argument that we will use repeatedly below.
\begin{lemma}
\label{lemm:ma}
Let $A$ be a commutative ring, $M$ an $A$-module, $a_u \in A_{\Qq}$ and $m_u \in M_{\Qq}$ for $u=0,\cdots,i(p-1)$ satisfying
\[
v_p(a_u) \geq -\Big[\frac{u}{p-1}\Big] \quad \text{ and } \quad v_p(m_u) \geq -\Big[\frac{u}{p-1}\Big].
\]
Then we have an equality in $M_{\Zz_{(p)} \subset \Qq} \otimes \Zz/p$,
\[
\sum_{u=0}^{i(p-1)} \big(p^i \cdot a_u \cdot m_{i(p-1)-u}\big) = \sum_{j=0}^{i} \big(p^j \cdot a_{j(p-1)}\big) \cdot \big(p^{i-j} \cdot m_{(i-j)(p-1)}\big).
\]
\end{lemma}
\begin{proof}
We compute, for all $u$, using \eqref{eq:vpmod},
\begin{align*}
v_p(a_u \cdot m_{i(p-1)-u})&\geq v_p(a_u) + v_p(m_{i(p-1)-u})\\
&\geq -\Big[\frac{u}{p-1}\Big] -\Big[\frac{i(p-1)-u}{p-1}\Big]\\
&=-\Big[\frac{u}{p-1}\Big] -i - \Big[\frac{-u}{p-1}\Big].
\end{align*}
This integer is always $\geq -i$, hence the elements delimited by parenthesis in the formula of the statement belong to the submodule $M_{\Zz_{(p)} \subset \Qq} \subset M_{\Qq}$. Moreover when $u$ is not divisible by $p-1$, this integer is $>-i$. This gives the requested formula modulo $p$.
\end{proof}

\begin{lemma}
\label{lemm:toddvp}
We have, for all integers $n$, and all $x \in \K(X)$
\[
v_p(\Td^n(x)) \geq -\Big[\frac{n}{p-1}\Big].
\]
\end{lemma}
\begin{proof}
Let $X$ be a variety. Then the set of elements $1+y$ with
\[
y \in \bigoplus_{n > 0} p^{-[n/(p-1)]}\cdot \CHc^n(X)_{\Zz_{(p)} \subset \Qq}
\]
is a subgroup $B(X)$ of $(\CHc(X)_{\Qq})^{\times}$. If $f\colon Y\to X$ is a morphism of varieties, then $f^*B(X) \subset B(Y)$. Thus we can define a presheaf of abelian groups $C$ on $\Var/k$ by
\[
C(-)=(\CHc(-)_{\Qq})^{\times}/ B(-).
\]
We have a morphism of presheaves of abelian groups on $\Var/k$
\[
\nu \colon \K(-) \to C(-), \quad x \mapsto \Td(-x).
\]
When $L$ is a line bundle, we have by definition of the Todd class,
\[
\nu[L]=\sum_{n \geq 0} \frac{(-c_1(L))^n}{(n+1)!}.
\]
Then $\nu$ vanishes on the classes of line bundles, because of \autoref{lemm:legendre} below. The splitting principle allows us to conclude.
\end{proof}

\begin{proposition}
\label{lemm:toddp}
Let $j \geq 0$ be an integer, $X$ a variety, and $x \in \K(X)$. Then the elements $p^j\cdot \Td^{j(p-1)}(-x)$ and $r^{(p)}_{j}(x)$ have same image in
\[
\CHc^{j(p-1)}(X)_{\Zz_{(p)} \subset \Qq} \otimes \Zz/p.
\]
\end{proposition}
\begin{proof}
We consider the presheaf of commutative rings on $\Var/k$
\[
\widetilde{\CHc}(-)=\CHc(-)_{\Zz_{(p)}\subset \Qq} \otimes \Zz/p.
\]
The characteristic class $r^{(p)}$ induces a morphism of presheaves of abelian groups on $\Var/k$
\[
\widetilde{r^{(p)}} \colon \K(-) \to \widetilde{\CHc}(-)^{\times}.
\]

By \autoref{lemm:toddvp}, we can consider the association
\[
\rho \colon \K(-) \to \widetilde{\CHc}(-), \quad x\mapsto \sum_{j \geq 0} p^j\cdot \Td^{j(p-1)}(-x).
\]
It is clearly compatible with arbitrary pull-backs. It sends sums to products because the total Todd class has the same property, and because of \autoref{lemm:ma} and \autoref{lemm:toddvp}. 

Let $L$ be a line bundle. Then by definition of the Todd class, we have
\[
\rho[L]=\sum_{j\geq 0} \frac{p^j}{(j(p-1)+1)!} \cdot (-c_1(L))^{j(p-1)}.
\]

Using \autoref{lemm:legendre} below, we see that the $j$-th coefficient is an element of $\Zz_{(p)}$. Moreover we see that it belongs to $p\Zz_{(p)}$ if and only if $j(p-1)+1$ is not a power of $p$. Otherwise, when say $j(p-1)+1=p^i$, we see using \autoref{lemm:pp} below, that the $j$-th coefficient is $(-1)^i$ modulo $p\Zz_{(p)}$. In view of \eqref{eq:r}, the characteristic class $\rho$ coincides with $\widetilde{r^{(p)}}$ on the classes of line bundles, and the claim follows from the splitting principle.
\end{proof}

\begin{lemma}[{\cite[p.10]{Legendre}}]
\label{lemm:legendre}
Let $p$ be a prime number, and $n$ a positive integer. Let $s$ the sum of the digits appearing in the expression of $n$ in base $p$. Then
\[
v_p(n!)=\frac{n-s}{p-1}.
\]
\end{lemma}

\begin{lemma}
\label{lemm:pp}
Let $i\geq 0$ be an integer, and $p$ a prime number. We have
\[
(p^i)! \cdot p^{-\frac{p^i-1}{p-1}}=(-1)^i \mod p.
\]
\end{lemma}
\begin{proof}
We proceed by induction on $i$. For $i=0$ this is immediate. Let $i\geq 1$ and consider the set $U_i=\{1,2,\cdots,p^i-1,p^i\}$. The product of the all the integers sitting between two consecutive multiples of $p$ is congruent to $-1$ modulo $p$, by Wilson's Theorem. Thus the product of all elements of $U_i$ which are non-divisible by $p$ is congruent modulo $p$ to $(-1)^{p^{i-1}}=-1$. The product of the $p^{i-1}$ elements of $U_i$ which are divisible by $p$ is $p^{p^{i-1}}\cdot (p^{i-1}) !$. Hence we have
\[
(p^i)! \cdot p^{-\frac{p^i-1}{p-1}}=(-1)\cdot p^{p^{i-1}}\cdot (p^{i-1})! \cdot p^{-\frac{p^i-1}{p-1}} = -p^{-\frac{p^{i-1}-1}{p-1}} \cdot (p^{i-1})!.
\]
By the induction hypothesis, this is $(-1)^i$, as requested.
\end{proof}

\section{Properties}
Let $p$ be a prime number. In this section, we assume that integrality holds in codimension $\leq m(p-1)$, and describe some functorial properties of the operations $T_0,\cdots,T_m$ constructed in \autoref{sect:construction}.

\begin{proposition}
\label{prop:compatpull}
Let $f \colon Y \to X$ be a local complete intersection morphism. Let $i$ an integer such that $0 \leq i\leq m$. Then we have
\[
T_i \circ f^*=\sum_{j=0}^i r^{(p)}_{j}(-\Tan_f) \cdot f^* \circ T_{i-j}.
\]
\end{proposition}
\begin{proof}
We have, by \cite[Theorem~18.3, (4)]{Ful-In-98},
\[
\ch_{n-i(p-1)} \circ f^*=\sum_{u=0}^{i(p-1)} \Td^u(\Tan_f) \cdot f^* \circ \ch_{n-i(p-1)+u}
\]
In view of \autoref{lemm:toddvp}, \autoref{lemm:ma}, and \autoref{th:mainvp}, we have, in the group $\CHZpn{n-i(p-1)}{Y} \otimes \Zz/p$,
\[
p^i \cdot \ch_{n-i(p-1)} \circ f^*=\sum_{j=0}^{i} \big(p^j \cdot \Td^{j(p-1)}(\Tan_f)\big) \cdot f^* \circ \big(p^{i-j}\cdot \ch_{n-(i-j)(p-1)}\big).
\]
Since the map \eqref{eq:chgrk} is compatible with $f^*$, \autoref{lemm:toddp} allows us to conclude.
\end{proof}

Applying \autoref{prop:compatpull} to the structural morphism of a local complete intersection variety, we obtain the formula computing the effect of the operation $T_i$ on regular cycles.
\begin{corollary}
\label{cor:reg}
Let $X$ be a local complete intersection variety (e.g. a regular variety). Let $i$ an integer such that $0\leq i\leq m$. Then
\[
T_i[X]=r_i^{(p)}(-\Tan_X)\cdot[X].
\]
\end{corollary}

\begin{proposition}
\label{prop:ext}
Let $X$ and $Y$ be two varieties over the same field. Then for $x \in \CH_d(X)$, and $y\in \CH_e(Y)$, and $0\leq i\leq m$, we have
\[
T_i(x\times y)=\sum_{j=0}^i T_j(x) \times T_{i-j}(y).
\]
\end{proposition}
\begin{proof}
Let $a \in \G(X)_{(d)}$ (resp. $b \in \G(X)_{(e)}$) be an element such that $\varphi_X(x)=a \mod  \G(X)_{(d-1)}$ (resp. $\varphi_Y(y)=b \mod  \G(Y)_{(e-1)}$), where $\varphi_{-}$ is the map \eqref{eq:chgrk}. We have, by \cite[Example~18.3.1]{Ful-In-98}, (and \cite[Theorem~18.3, (5)]{Ful-In-98})
\[
\ch_{d+e-i(p-1)}(a \times b)=\sum_{u=0}^{i(p-1)} \ch_{d-u}(a) \times \ch_{e+u-i(p-1)}(b).
\]
By \autoref{lemm:ma}, this gives in $\CHZpn{d+e-i(p-1)}{X\times Y} \otimes \Zz/p$,
\[
p^i\cdot \ch_{d+e-i(p-1)}(a \times b)=\sum_{j=0}^{i} \big(p^j\cdot \ch_{d-j(p-1)}(a)\big) \times \big(p^{i-j}\cdot \ch_{e-(i-j)(p-1)}(b)\big).
\]
We can conclude using the fact that $\varphi$ is compatible with external products.
\end{proof}

Assume that $X$ is a local complete intersection variety. Define, for $i\leq m$, the \emph{$i$-th cohomological operation} $T^i$ by the formula 
\[
T^i=\sum_{j=0}^i r_{j}^{(p)}(\Tan_X) \cdot T_{i-j} \quad \colon \CH_\bullet(X) \to \CHZn{\bullet-i(p-1)}{X}\otimes \Zz/p.
\]

\begin{proposition}
\label{prop:coh}
Let $X$ and $Y$ be local complete intersection varieties over a common field, and $i$ an integer such that $0<i\leq m$.  
\begin{enumerate}[(i)]
\item \label{coh:ext} For $x \in \CH(X), y\in\CH(Y)$, we have 
\[
T^i(x \times y)=\sum_{j=0}^i T^j(x) \times T^{i-j}(y).
\]

\item \label{coh:prod} When $X$ is smooth, $x,y \in \CH(X)$, we have 
\[
T^i(x\cdot y)=\sum_{j=0}^i T^j(x)\cdot T^{i-j}(y).
\]

\item \label{coh:id} The operation $T^0$ is the identity.

\item \label{coh:i} $T^i[X]=0$. 

\item \label{coh:one} Assume that $X$ is smooth, and let $x \in \CH^1(X)$. Then $T^i(x)=-x^p$ if $i=1$, and $T^i(x)=0$ otherwise.

\item \label{coh:pull} If $f \colon Y \to X$ is a local complete intersection morphism, then $T^i \circ f^*=f^* \circ T^i$.

\item \label{coh:pushforward} Let $f\colon Y \to X$ be a projective morphism. Let $\Tan_f=\Tan_Y-f^*\Tan_X \in \K(Y)$. Then
\[
T^i \circ f_*=f_* \circ \sum_{j=0}^i r^{(p)}_j(-\Tan_f) \cdot T^{i-j}.
\]
\end{enumerate}
\end{proposition}
\begin{proof}
Statement \eqref{coh:ext} follows from \autoref{prop:ext}, statement \eqref{coh:pull} from \autoref{prop:compatpull}, and \eqref{coh:pushforward} from the fact that the operations $T_i$ are compatible with projective push-forwards. Then \eqref{coh:prod} follows from \eqref{coh:ext} and \eqref{coh:pull}, statement \eqref{coh:i} from \eqref{coh:pull} and the obvious computation for $X=\Spec(k)$. Statement \eqref{coh:id} is a consequence of the definition of $T_0$. 

We now prove \eqref{coh:one}. We can assume that $x=[Z]$ for some integral closed subvariety $Z$ of $X$. Since $X$ is locally factorial, $Z$ is a locally principal divisor. In particular the closed embedding $f\colon Z \hookrightarrow X$ is regular, and $Z$ is a local complete intersection variety. Then, using \eqref{coh:i} and \eqref{coh:pushforward}, where $T_f=-f^*[\Oc_X(Z)]$, we obtain for $i>0$,
\begin{align*}
T^i(x)=T^i\circ f_*[Z]&=f_* \circ T^i[Z] - f_* (c_1(f^*[\Oc_X(Z)])^{p-1} \cdot T^{i-1}[Z])\\
&=-c_1(\Oc_X(Z))^{p-1}\cdot f_*\circ T^{i-1}[Z].
\end{align*}
This is zero when $i>1$ by \eqref{coh:i}. When $i=1$, this is, by \eqref{coh:id},
\[
-c_1(\Oc_X(Z))^{p-1}\cdot [Z]=-c_1(\Oc_X(Z))^p=-x^p,
\]
as requested.
\end{proof}

Finally, we mention
\begin{proposition}
\label{prop:fieldext}
Let $L/k$ be a field extension. Then, for all varieties $X$, all $x \in \CH(X)$, and $i$ such that $0\leq i\leq m$, we have
\[
T_i(x_L)=T_i(x)_L.
\]
\end{proposition}
\begin{proof}
This follows from \cite[(2) and Theorem~3.1, (d)]{firstsq}.
\end{proof}

\section{Relation with Steenrod operations}
\label{sect:relations}
In this section we compare the operations $T_i$ with the Steenrod operations constructed by other methods. We therefore make the assumption that the characteristic of the base field is different from $p$.\\

Consider the $i$-th homological Steenrod operation of, say, \cite{Bro-St-03}
\[
S_i \colon \CH_\bullet(-) \otimes \Zz/p \to \CH_{\bullet-i(p-1)}(-) \otimes \Zz/p,
\]
and write $S=\sum_i S_i$. Since $S_0=\id$, we have $S=\id+f$, where for every variety $X$ the induced endomorphism $f^X$ of $\CH(X) \otimes \Zz/p$ is nilpotent. We write $\Te^X$ for the left inverse of $S^X$; this defines a natural transformation of functors $\Var/k \to \Ab$ with components
\[
\Te_i \colon \CH_\bullet(-) \otimes \Zz/p \to \CH_{\bullet-i(p-1)}(-) \otimes \Zz/p.
\]
We have  $\Te_0=\id$ and, for $i>0$, one can express $\Te_i$ inductively as
\begin{equation}
\label{eq:induct}
\Te_i=-\sum_{j=1}^i \Te_{i-j} \circ S_j.
\end{equation}

\begin{proposition}
\label{prop:st}
We have for $i\leq p$
\[
\Te_i=(-1)^i \cdot S_i.
\]
\end{proposition}
\begin{proposition}
\label{prop:tr}
Assume that $X$ is a smooth variety. Then
\[
T'[X]=r^{(p)}(-\Tan_X)[X].
\]
\end{proposition}
\begin{proof}[Proof of \autoref{prop:st} and \autoref{prop:tr}]
Let $Z$ be a variety. It is proven in \cite[Lemma~5.2 and Proposition~3.1]{Mer-St-03} that the \emph{homological} Steenrod operations $S_\bullet$ induce an action of the Steenrod $\Zz/p$-algebra $\mathcal{S}$ on $\CH(Z) \otimes \Zz/p$. The operation $S_j$ acts as the element $\mathcal{P}^j\in \mathcal{S}$ (in case $p=2$ we write $\mathcal{P}^j$ for $\Sq^{2j}$) . Using \eqref{eq:induct}, we see that the action of $\Te_j$ on $\CH(Z) \otimes \Zz/p$ coincides with the action of the conjugate $c(\mathcal{P}^j)$ of $\mathcal{P}^j$, in the sense of \cite[\S 7]{Milnor-Steenrod}. Let $b_i$ for $i \in \mathbb{N}-\{0\}$ be indeterminates. If $R=(r_1,r_2,\cdots)$ is a sequence of non-negative integers, we write $b^R$ for the product ${b_1}^{r_1}{b_2}^{r_2}\cdots$. By \cite[Corollary 6]{Milnor-Steenrod}, the morphism
\[
\mathcal{S}[ [b_1,b_2,\cdots] ] \to \mathcal{S}[[x]], \quad \quad  b_i \mapsto (-1)^{\frac{p^i-1}{p-1}}\cdot x^{p^i-1}=(-1)^i \cdot x^{p^i-1} \mod p
\]
sends the total Steenrod operation $\sum_R \mathcal{P}^R \cdot b^R$ to $\sum_j c(\mathcal{P}^j)\cdot x^{j(p-1)}$.

We obtain \autoref{prop:st} by noticing that for $i\leq p$ the element $\mathcal{P}^i$ is the only basis element of type $\mathcal{P}^R$ of the same dimension as $c(\mathcal{P}^i)$.

Moreover by \cite[Proposition~5.3, (2)]{Mer-St-03}, we obtain
\[
T'[X]=d(-\Tan_X)[X],
\]
where $d$ is the total Chern class $\sum_R c_R \cdot b^R$ evaluated at $b_i=(-1)^i$. Inspection of the value of $d$ on line bundles reveals that $d=r^{(p)}$, establishing \autoref{prop:tr}.
\end{proof}

\begin{proposition}
\label{prop:smooth}
Assume that the base field is of characteristic different from $p$. Let $i$ be any integer $\geq 0$. Then we have a commutative diagram
\[ \xymatrix{
\CH(-) \otimes \Zz/p\ar[r]^{T_i'} \ar[dr]_{T_i} &  \CH(-) \otimes \Zz/p\ar[d] \\ 
& \CHZ{-} \otimes \Zz/p
}\]
\end{proposition}
\begin{proof}
By \autoref{lemm:perfectclosure} below and \autoref{prop:fieldext}, we can assume that the base field is perfect.

Since both operations are compatible with proper push-forwards, it will be enough to compare their values in $\CHZn{d-i(p-1)}{X} \otimes \Zz/p$ on the class $[X]$, when $X$ is an integral variety of dimension $d$. Using Gabber's theorem, we find a projective morphism $f \colon Z \to X$, with $Z$ smooth of dimension $d$, and $f$ dominant with $[k(Z):k(X)]=\lambda$ an integer prime to $p$. Let $\nu \in (\Zz/p)^{\times}$ be the inverse modulo $p$ of $\lambda$. Then, using \autoref{cor:reg},
\[
T_i[X]=\nu\cdot T_i(\lambda \cdot [X])=\nu \cdot T_i \circ f_*[Z]=\nu \cdot f_* \circ T_i[Z]=\nu \cdot f_*  (r_i^{(p)}(-\Tan_Z)\cdot[Z]).
\]
A similar computation, and \autoref{prop:tr}, show that $T_i'[X]$ has the same value, as requested.
\end{proof}
\begin{lemma}
\label{lemm:perfectclosure}
Let $X$ be a variety over a field $k$ of characteristic $\neq p$. Let $F/k$ be a perfect closure. Then scalars extension induces an injective morphism 
\[
\CHZ{X}\otimes \Zz/p \to \CHZ{X_F}\otimes \Zz/p.
\]
\end{lemma}
\begin{proof}
Let $L/k$ be a finite field extension. A transfer argument shows that $\CHQ{X} \to \CHQ{X_L}$ is injective, hence so is $\CHZ{X} \to \CHZ{X_L}$, and therefore so is
\[
\CHZ{X}\otimes \Zz_{(p)} \to \CHZ{X_L}\otimes \Zz_{(p)}.
\]
When the degree of the extension $L/k$ is prime to $p$, the transfer induces a splitting of this monomorphism, so that
\[
\CHZ{X}\otimes \Zz/p \to \CHZ{X_L}\otimes \Zz/p
\]
is injective.

By \cite[Theorem~5.20 and Proposition~5.26]{Sri-96}, the group $\CH(X_F)$ is the direct limit of the groups $\CH(X_L)$ with $L/k$ finite subextension of $F/k$. Such an extension $L/k$ is of degree prime to $p$ because it is a purely inseparable extension of a field of characteristic $\neq p$. Using the fact that the functor $A \mapsto A_{\Zz \subset \Qq} \otimes \Zz/p$ commutes with any exact functor commuting with tensor products, we see that the natural map 
\[
\colim{\substack{k\subset L\subset F\\ [L:k]<\infty}} (\CHZ{X_L}\otimes \Zz/p) \to \CHZ{X_F}\otimes \Zz/p
\]
is an isomorphism. The lemma follows.
\end{proof}

\begin{remark}
\label{rem:K}
\autoref{prop:smooth} actually gives some new informations about Steenrod operations modulo $p$ in characteristic not $p$. Denote by $\mu\colon \CH(-) \otimes \Zz/p \to \CHZ{-} \otimes \Zz/p$ the natural quotient map. Let $X$ be a variety, and
\begin{equation}
\label{eq:K}
\mathcal{K}(X)=\ker(\CH(X) \to \gr \G(X)).
\end{equation}
Then $\mathcal{K}(X)$ consists of torsion elements. Let $\mathcal{K}_p(X)$ be its image in $\CH(X) \otimes \Zz/p$. One can ask the following question: ``what is the set $\mathcal{E}$ of operations $Q$ on Chow groups modulo $p$ such that $\mu \circ Q$ vanishes on $\mathcal{K}_p(X)$, for all $X$?''

By construction, we have $T_i(\mathcal{K}_p(X))=0$, hence a consequence of \autoref{prop:smooth} is that $T_i' \in \mathcal{E}$, for all $i$. On the other hand when $p=2$ and $j>1$, then the cohomological Steenrod square $S^{2^j-1}$ does not belong to $\mathcal{E}$ (see \cite[Remark~8.4]{duality}).

This can also be used to obtain informations on $\gr \G(X)$. For example, a direct consequence of \autoref{prop:smooth} is that $\gr_2 \G(D)$ is contains non-zero $3$-torsion, where $D$ is the variety of \cite{Semenov-compactification} (over a field of characteristic zero).
\end{remark}

\begin{remark}
One can use the operations $T_i$ to recover some of the results of \cite{reduced}. Assume that $X$ is a variety and $L/k$ an extension of the base field such that $\CH(X_L)$ is torsion-free. Then $T_i$ is an endomorphism of $\CH(X_L)\otimes \Zz/p$. We set $S_0=\id$, and define inductively, for $n>0$ and $x \in \CH(X_L)\otimes \Zz/p$,
\[
S_n(x)=-\sum_{i=1}^n T_i \circ S_{n-i}(x).
\]
Let $\overline{\Ch}(X)=\im(\CH(X)\otimes \Zz/p \to \CH(X_L)\otimes \Zz/p)$ be the reduced Chow group modulo $p$. Since the operations $T_i$ commute with extension of the base field, they send the subgroup $\overline{\Ch}(X)\subset  \CH(X_L)\otimes \Zz/p$ to itself. Therefore the association $x \mapsto S_n(x)$ induces an operation of this subgroup, the $n$-th \emph{reduced Steenrod operation}. 
\end{remark}

\part{Degree formula and applications}
We now describe some consequences of the results obtained above. In order to make unconditional statements, we restrict ourselves to applications of the fact that $p$-integrality holds in codimension $\leq p(p-1)-1$ (\autoref{th:mainvp}).\\

We fix prime number $p$. When $X$ is a projective variety, we let $n_X$ be the positive generator of the image of the degree map $\CH(X) \to \Zz$. We also set
\[
\nx{X}=p^{v_p(n_X)}.
\]
Let $x\colon X \to \Spec(k)$ be the structural morphism of $X$. We denote by $\chi(\Oc_X)$ the integer corresponding to $x_*[\Oc_X]$ under the identification $\G(\Spec(k))=\Zz$. In other words
\[
\chi(\Oc_X)=\sum_{j \geq 0} (-1)^j \dim_k H^j(X,\Oc_X).
\]
Note that $\ch[\Oc_{\Spec(k)}]=[\Spec(k)]$, so that 
\[
\chi(\Oc_X)=\deg \circ \ch_0[\Oc_X].
\]

\section{Index and Euler characteristic}
\label{sect:degree}
We begin by giving a restriction on the possible values of the index $n_X$ of a projective variety $X$, in terms of two integers which do not change under field extensions: the dimension of $X$ and the Euler characteristic $\chi(\Oc_X)$.

\begin{proposition}
\label{prop:majornx}
Let $X$ be a projective variety of dimension $< p(p-1)$. Then
\[
v_p(n_X) \leq \Big[ \frac{\dim X}{p-1} \Big] + v_p(\chi(\Oc_X)).
\]
\end{proposition}
\begin{proof}
We have, by \autoref{th:mainvp}, 
\[
v_p(\ch_0[\Oc_X]) \geq -\Big[ \frac{\dim X}{p-1} \Big],
\]
which means that there is an integral zero-cycle $x\in \CH_0(X)$, and an integer $\lambda$ prime to $p$ such that
\[
x=\lambda \cdot p^{[\dim X/(p-1)]} \cdot \ch_0[\Oc_X].
\]
Taking degrees, we obtain
\[
n_X | \lambda \cdot p^{[\dim X/(p-1)]} \cdot \chi(\Oc_X),
\]
whence the result.
\end{proof}

Let $X$ be an integral projective variety of dimension $i(p-1)$, with $i>0$. We define in $\Zz/\nx{X}$
\[
t_p(X)=p^{i-1} \cdot \chi(\Oc_X) \mod \nx{X}.
\]
Then \autoref{prop:majornx} amounts to saying that we have, when $i< p$,
\[
p \cdot t_p(X)=0 \mod \nx{X}.
\]
The same conclusion holds for any $i>0$ when $X$ is a local complete intersection variety by \autoref{prop:smooth2}, or if $p$-integrality holds in codimension $\leq i(p-1)$.

\begin{remark}
By \cite[Corollary~1]{Ros-On-08}, the class $t_p(X)$ vanishes when the base field is perfect of characteristic $p$. 
\end{remark}

\begin{remark}
\label{rem:Rost}
Assume that $i < p$. Let $u_p(X) \in \CH_0(X)$ be a cycle whose class modulo $p$ and $p$-primary torsion is $T_i[X]$. Then its image $\overline{u_p}(X)$ in $\CHZpn{0}{X}$ satisfies
\[
\overline{u_p}(X)=p^i\cdot \ch_0[\Oc_X] \mod p \CHZpn{0}{X}.
\]
Hence
\[
\deg \overline{u_p}(X)=p^i \cdot \chi(\Oc_X) \mod p n_X \Zz_{(p)}
\]
and
\[
\frac{\deg u_p(X)}{p}=p^{i-1} \cdot \chi(\Oc_X) \mod n_X(p).
\]
Thus the class $t_p(X)$ corresponds to the Rost numbers for a certain combination of Steenrod operations (see \cite[p.20]{Mer-St-03}) under the map $\Zz/n_X \to \Zz/n_X(p)$. 
\end{remark}

When $f\colon Y \to X$ is a morphism between varieties of the same dimension, then the integer $\deg f$ is defined as the degree of the extension of function fields when $f$ is dominant, and as zero otherwise.
\begin{theorem}[Degree formula]
\label{th:degf}
Let $f \colon Y \to X$ be a morphism of integral projective varieties of the same dimension $i(p-1)$ with $0<i\leq p$. Then
\[
t_p(Y)=\deg f \cdot t_p(X) \mod \nx{X}.
\]
\end{theorem}
\begin{proof}
We have in $\G(X)$
\begin{equation}
\label{eq:deg}
f_*[\Oc_Y]=\deg f \cdot [\Oc_X] + \delta
\end{equation}
with $\delta \in \G(X)_{(i(p-1)-1)}$. This element $\delta$ can be written as a linear combination of elements of type $[\Oc_Z]$, for $Z$ a closed subvariety of $X$, with $\dim Z\leq i(p-1)-1$. Note that $n_X(p)|n_Z(p)$ for such $Z$. The result is obtained by projecting \eqref{eq:deg} to $\G(\Spec(k))$, and using \autoref{prop:majornx}.
\end{proof}

Thus the class $t_p(-)$ is defined for varieties of dimension $i(p-1)$ when $i$ is arbitrary, satisfies the degree formula when $i\leq p$, and is of $p$-torsion when $i \leq p-1$.

\begin{remark}
Let $X$ be a projective variety of dimension $i(p-1)$, with $i \leq p $. Assume that there is a morphism $f \colon X \to Y$, for some projective variety $Y$, such that $f_*[X] \in \mathcal{K}(Y)$ (see \eqref{eq:K}). Then $t_p(X)=0 \mod n_Y(p)$.
\end{remark}

\section{Correspondences modulo $p$}
Let $X$ and $Y$ be varieties, with $Y$ projective and $X$ integral. A correspondence $\alpha \colon X \leadsto Y$ is an element $\alpha \in \CH(X \times Y)$. The \emph{multiplicity} of $\alpha$ is the image of $\alpha$ under the composite 
\[
\CH(X \times Y) \to \CH(k(X) \times Y) \xrightarrow{\deg} \CH(k(X))=\Zz.
\]

\subsection{Strong $p$-incompressibility}
Let $p$ be a prime number. Following \cite[p.150]{Kar-can}, we call an integral projective variety $X$ \emph{strongly $p$-incompressible} if the existence of a correspondence $X \leadsto Y$ of multiplicity prime to $p$, with $Y$ an integral projective variety of dimension $\leq \dim X$ and such that $v_p(n_Y) \geq v_p(n_X)$, implies that $\dim Y = \dim X$, and that there exists a correspondence $Y \leadsto X$ of multiplicity prime to $p$.

\begin{proposition}
\label{prop:m}
Let $p$ be a prime number, and $i$ an integer such that $0 < i \leq p$. Let $X$ be an integral projective variety with $v_p(n_X)\geq i$. Assume that $\chi(\Oc_X)$ is prime to $p$. Then $\dim X\geq i(p-1)$. In case of equality, $X$ is strongly $p$-incompressible. 
\end{proposition}
\begin{proof}
Assume that $\dim X<i(p-1)$. Then in particular $\dim X <p(p-1)$, and \autoref{prop:majornx} gives 
\[
i\leq v_p(n_X) \leq \Big[ \frac{\dim X}{p-1} \Big] + v_p(\chi(\Oc_X)) < i,
\]
a contradiction. Hence $\dim X \geq i(p-1)$.

Now we assume that $\dim X=i(p-1)$. We proceed as in the proof of \cite[Theorem~7.2]{Mer-St-03} --- but as we include the case $i=p$, we cannot use the relation $p\cdot t_p(X)=0$.

Let $Y$ be an integral projective variety of dimension $\leq \dim X$, with a correspondence $X\leadsto Y$ of multiplicity prime to $p$, and such that $v_p(n_Y)\geq v_p(n_X)$.  We can assume that there is an integral variety $Z$ with $\dim Z=\dim X$, and two projective morphisms $f\colon Z \to X$ and $g \colon Z \to Y \hookrightarrow Y'$, such that $\deg f$ is prime to $p$. Here $Y'$ is the product of $Y$ with a projective space of appropriate dimension: in particular $\dim Y'=\dim X$ and $n_{Y'}=n_Y$.  Since by \autoref{th:degf} we have
\[
\deg f \cdot t_p(X)=t_p(Z) \mod n_X(p),
\]
we see that $t_p(Z)$ is non-zero modulo $n_X(p)$, hence modulo $n_Y(p)$. Now by the same theorem we have
\[
\deg g \cdot t_p(Y')=t_p(Z) \mod n_Y(p),
\]
hence $\deg g$ is non-zero and $g$ is dominant, so that $Y=Y'$ and $\dim Y=\dim X$. Then $Z$ gives a correspondence $Y \leadsto X$ of multiplicity prime to $p$.

For further reference, let us mention that we also obtain that $t_p(Y)$ is non-zero modulo $n_Y(p)$.
\end{proof}

Let $X$ be the Severi-Brauer variety of a central division algebra of degree $p^n$, with $p$ a prime number. Then $X$ is strongly $p$-incompressible when the characteristic of the base field is not $p$ (\cite[Example~7.2, Theorem~7.2]{Mer-St-03}). According to \cite[Example~2.3]{Kar-can}, it is unknown whether $X$ is strongly $p$-incompressible when the characteristic is $p$. We can treat the case $n=1$ with \autoref{prop:m}:
\begin{corollary}
The Severi-Brauer variety of a central division algebra of prime degree $p$ is strongly $p$-incompressible.
\end{corollary}

\subsection{Index reduction formula}
We mentioned in \cite[\S 4.1]{Euler} that the index reduction formula of \cite[Corollary (B)]{Zai-09} can be generalized to base fields of arbitrary characteristic. The results obtained here allow us to additionally remove all assumptions of regularity:
\begin{proposition}
Let $X$ and $Y$ be integral projective varieties, of dimensions $\leq i(p-1)$ for some $i\leq p$. Assume that $v_p(n_Y)\geq v_p(n_X) \geq i$, that $p | \chi(\Oc_Y)$, and that $\chi(\Oc_X)$ is prime to $p$. Then there is no correspondence $X \leadsto Y$ of multiplicity prime to $p$.
\end{proposition}
\begin{proof}
The first part of \autoref{prop:m} implies that $\dim X=i(p-1)$, hence $\dim Y \leq \dim X$. If such a correspondence existed, the last sentence of the proof of \autoref{prop:m} would imply that $\chi(\Oc_Y)$ is prime to $p$, a contradiction.
\end{proof}

\subsection{Rost's correspondences}
The operations $T_i$ can be used to extend the case $n=1$ in \cite[Theorem~9.1]{Ros-On-06} to arbitrary fields.

\begin{theorem}
Let $X$ be a smooth projective variety of dimension $p-1$. Assume that $X$ has a special correspondence in the sense of \cite[Definition~5.1]{Ros-On-06}, with $b=1$. Assume moreover that $X$ has no zero-cycle of degree prime to $p$. Then 
\[
\deg (r^{(p)}(-\Tan_X)\cdot [X])\neq 0 \mod p^2.
\]
By \autoref{rem:Rost}, this means that the integer $\chi(\Oc_X)$ is prime to $p$.
\end{theorem}
\begin{proof}
We only sketch how \cite[Proof of Theorem~9.1]{Ros-On-06} can be adapted. This proof only involve cycles of dimension $\leq 2p-3$, in particular the only Steenrod operation used is the first one. Moreover since the degree map factors through the Chow group modulo torsion, there is no loss in using the operation with values in this group. All the properties needed are contained in \autoref{prop:coh}.
\end{proof}

\bookmarksetup{startatroot}
\addtocontents{toc}{\bigskip}
\appendix
\section{Schemes over a regular base}
\label{sect:base}
Let $S$ be a regular connected scheme of finite Krull dimension. Then, as explained in \cite[\S 20.1]{Ful-In-98}, most results of \cite{Ful-In-98} for varieties over fields extend to arbitrary schemes of finite type over $S$. There are only two difficulties : 

--- the lack of external product in general.

--- one should use relative dimension $\dim_S$ instead of the usual dimension.

In particular the homological Chern character can be constructed in this generality. The purpose of this appendix is to explain how \autoref{th:mainvp} (and thus the construction of the operations $T_i$, $i=0,\cdots,p-1$) can be extended to this setting.\\

Let $\mathcal{V}_S$ be the category of schemes quasi-projective over $S$. Homological Adams operations are constructed in \cite{Sou-Op-85} for objects of $\mathcal{V}_S$. The results of \autoref{sect:adamschern} extend without difficulty to the category $\mathcal{V}_S$.

However, both proofs of \autoref{prop:graded} used the fact that $X$ contains a non-empty open regular subscheme. We now give a proof that works for arbitrary objects of $\mathcal{V}_S$, using the notations introduced in the beginning of \autoref{sect:adamschern}.

\begin{proposition}
\label{prop:graded2}
Let $l \in \Zz - \{0\}$, $X \in \mathcal{V}_S$, and $x \in \G(X)_{(d)}$. Then
\[
\psi_l(x)-l^{-d} \cdot x \in \G(X)_{(d-1)} \otimes \Zz[1/l].
\]
\end{proposition}
\begin{proof}
As before, one can assume that $X$ is integral and $d=\dim_S X$, and we have the exact sequence
\[
0 \to \G(X)_{(d-1)} \otimes \Zz[1/l] \to \G(X)_{(d)}\otimes \Zz[1/l] \xrightarrow{u^*} \G(\Spec(k(X)))\otimes \Zz[1/l] \to 0.
\]
Choose a closed embedding $X \hookrightarrow W$, with $W$ connected and smooth over $S$. Let $x$ be the generic point of $X$, and write $W_x$ for $\Spec(\Oc_{W,x})$. Then the map $u^*$ in the sequence above may be identified with the map
\[
f^* \colon \K^X(W)\otimes \Zz[1/l] \to  \K^{\{x\}}(W_x) \otimes \Zz[1/l]
\]
and the composite $u^* \circ \psi_l$ corresponds to
\[
f^* \circ (\theta^l(-\Tan_{W/S}) \cdot \psi^l_X)=\theta^l(-\Tan_{W/S}) \cdot \psi^l_{\{x\}} \circ f^*.
\]
The action of $\theta^l(-\Tan_{W/S}) \in \K(W)$ on $\K^{\xx}(W_x)$ factors through $\K(W_x)$. Since $W_x$ is a local scheme, the rank map $\K(W_x) \to \mathbb{Z}$ is an isomorphism. The rank of $\theta^l(-\Tan_{W/S}))~$ is $l^{-\rank\Tan_{W/S}}$ (\cite[(2), p.6]{duality}). Since $W$ is flat over $S$, we have $\dim_S W=\trdeg(k(W)/k(S))=\rank \Tan_{W/S}$. Thus $\theta^l(-\Tan_{W/S}))$ acts on $\K^{\xx}(W_x)$ by multiplication by $l^{-\dim_S W}$. By \cite[A4') p.263, see also Proof of Proposition~5.3]{GS-87}, $\psi^l_{\xx}$ is multiplication by $l^{\codim(X,W)}$. This gives, using \cite[Lemma~20.1, (2)]{Ful-In-98},
\[
u^* \circ \psi_l=l^{\codim(X,W)-\dim_S W} \cdot u^*=l^{-\dim_S X} \cdot u^*.
\]
We can now conclude using the exact sequence above.
\end{proof}

As a consequence, \autoref{th:mainvp} immediately extends to the category $\mathcal{V}_S$. 

Finally let $X$ a scheme of finite type over $S$ (possibly non quasi-projective). It follows from \cite[Lemma~18.3, (1)]{Ful-In-98} that there is always an envelope $X' \to X$ over $S$, such that $X'\in \mathcal{V}_S$. Moreover we see, using \cite[Lemma~20.1, (3)]{Ful-In-98}, that $\G(X')_{(d)} \to \G(X)_{(d)}$ is surjective, for any $d$. Using this remark, \autoref{th:mainvp} is easily extended to the category of schemes of finite type over $S$.\\

{\bf Acknowledgements.} The support of EPSRC Responsive Mode grant EP/G032556/1 is gratefully acknowledged. I thank Alexander Vishik for his interest in this work, and for the useful suggestions that he made.

\bibliographystyle{alpha}

\begin{thebibliography}{A}

\bibitem[AA66]{Adams-Atiyah}
John~F. Adams and Michael~F. Atiyah.
\newblock {$K$}-theory and the {H}opf invariant.
\newblock {\em Quart. J. Math. Oxford Ser. (2)}, 17:31--38, 1966.

\bibitem[Ada61]{Ada-On-61}
John~F. Adams.
\newblock On {C}hern characters and the structure of the unitary group.
\newblock {\em Proc. Cambridge Philos. Soc.}, 57:189--199, 1961.

\bibitem[Ada62]{Adams-Vector}
John~F. Adams.
\newblock Vector fields on spheres.
\newblock {\em Ann. of Math. (2)}, 75:603--632, 1962.

\bibitem[Ati66]{Ati-Po-66}
Michael~F. Atiyah.
\newblock Power operations in {$K$}-theory.
\newblock {\em Quart. J. Math. Oxford Ser. (2)}, 17:165--193, 1966.

\bibitem[BFM75]{BFM}
Paul Baum, William Fulton, and Robert MacPherson.
\newblock Riemann-{R}och for singular varieties.
\newblock {\em Inst. Hautes \'Etudes Sci. Publ. Math.}, (45):101--145, 1975.

\bibitem[Boi08]{Boi-A-08}
Alex Boisvert.
\newblock A new definition of the {S}teenrod operations in algebraic geometry.
\newblock {\em preprint}, 2008.
\newblock \href{http://arxiv.org/abs/0805.1414}{\tt{arXiv:0805.1414}}.

\bibitem[Bro03]{Bro-St-03}
Patrick Brosnan.
\newblock Steenrod operations in {C}how theory.
\newblock {\em Trans. Amer. Math. Soc.}, 355(5):1869--1903 (electronic), 2003.

\bibitem[EKM08]{EKM}
Richard Elman, Nikita Karpenko, and Alexander Merkurjev.
\newblock {\em The algebraic and geometric theory of quadratic forms},
  volume~56 of {\em American Mathematical Society Colloquium Publications}.
\newblock American Mathematical Society, Providence, RI, 2008.

\bibitem[Ful98]{Ful-In-98}
William Fulton.
\newblock {\em Intersection theory}, volume~2 of {\em Ergebnisse der Mathematik
  und ihrer Grenzgebiete. 3. Folge. A Series of Modern Surveys in Mathematics
  [Results in Mathematics and Related Areas. 3rd Series. A Series of Modern
  Surveys in Mathematics]}.
\newblock Springer-Verlag, Berlin, second edition, 1998.

\bibitem[Gab]{Gabber}
Ofer Gabber.
\newblock Finiteness theorems in \'etale cohomology.
\newblock {\em Abstract of a talk given at the IHES on January 12, 2009}.
\newblock \url{http://www.ihes.fr/document?id=1713&id_attribute=48}.

\bibitem[Gil05]{Gil-K-05}
Henri Gillet.
\newblock {$K$}-theory and intersection theory.
\newblock In {\em Handbook of {$K$}-theory. {V}ol. 1, 2}, pages 235--293.
  Springer, Berlin, 2005.

\bibitem[GS87]{GS-87}
Henri Gillet and Christophe Soul{\'e}.
\newblock Intersection theory using {A}dams operations.
\newblock {\em Invent. Math.}, 90(2):243--277, 1987.

\bibitem[Haua]{Euler}
Olivier Haution.
\newblock Degree formula for the {E}uler characteristic.
\newblock {\em Proc. Amer. Math. Soc.}, to appear.
\newblock \href{http://arxiv.org/abs/1006.1482}{\tt{arXiv:1103.4076}}.

\bibitem[Haub]{firstsq}
Olivier Haution.
\newblock On the first {S}teenrod square for {C}how groups.
\newblock {\em Amer. J. Math.}, to appear.
\newblock \url{ http://www.math.uiuc.edu/K-theory/0945/}.

\bibitem[Hauc]{reduced}
Olivier Haution.
\newblock Reduced {S}teenrod operations and resolution of singularities.
\newblock {\em J. K-Theory}, to appear.
\newblock \href{http://arxiv.org/abs/1006.1480}{\tt{arXiv:1006.1480}}.

\bibitem[Hau10]{duality}
Olivier Haution.
\newblock Duality and the topological filtration.
\newblock {\em preprint}, 2010.
\newblock \href{http://arxiv.org/abs/1006.1482}{\tt{arXiv:1006.1482}}.

\bibitem[Ill]{Gabber-Illusie}
Luc Illusie.
\newblock On {G}abber's refined uniformization.
\newblock {\em Talks at the Univ. Tokyo, Jan. 17, 22, 31, Feb. 7, 2008}.
\newblock \url{http://www.math.u-psud.fr/~illusie/refined_uniformization3.pdf}.

\bibitem[Kar10]{Kar-can}
Nikita~A. Karpenko.
\newblock Canonical dimension.
\newblock {\em Proceedings of the ICM 2010}, Vol. II:146--161, 2010.

\bibitem[KR00]{Kurano-Roberts}
Kazuhiko Kurano and Paul~C. Roberts.
\newblock Adams operations, localized {C}hern characters, and the positivity of
  {D}utta multiplicity in characteristic {$0$}.
\newblock {\em Trans. Amer. Math. Soc.}, 352(7):3103--3116, 2000.

\bibitem[Leg08]{Legendre}
Adrien-Marie Legendre.
\newblock {\em Essai sur la th\'eorie des nombres}.
\newblock Courcier, Paris, second edition, 1808.

\bibitem[Lev07]{Lev-St-05}
Marc Levine.
\newblock Steenrod operations, degree formulas and algebraic cobordism.
\newblock {\em Pure Appl. Math. Q.}, 3(1):283--306, 2007.

\bibitem[Mer03]{Mer-St-03}
Alexander Merkurjev.
\newblock Steenrod operations and degree formulas.
\newblock {\em J. Reine Angew. Math.}, 565:13--26, 2003.

\bibitem[Mil58]{Milnor-Steenrod}
John Milnor.
\newblock The {S}teenrod algebra and its dual.
\newblock {\em Ann. of Math. (2)}, 67:150--171, 1958.

\bibitem[Pan04]{Pan-Ri-02}
Ivan Panin.
\newblock Riemann-{R}och theorems for oriented cohomology.
\newblock In {\em Axiomatic, enriched and motivic homotopy theory}, volume 131
  of {\em NATO Sci. Ser. II Math. Phys. Chem.}, pages 261--333. Kluwer Acad.
  Publ., Dordrecht, 2004.

\bibitem[Ros06]{Ros-On-06}
Markus Rost.
\newblock On the basic correspondence of a splitting variety.
\newblock {\em preprint}, 2006.
\newblock \url{http://www.math.uni-bielefeld.de/~rost/basic-corr.html}.

\bibitem[Ros08]{Ros-On-08}
Markus Rost.
\newblock On {F}robenius, ${K}$-theory, and characteristic numbers.
\newblock {\em preprint}, 2008.
\newblock \url{http://www.math.uni-bielefeld.de/~rost/frobenius.html}.

\bibitem[Sem08]{Semenov-compactification}
N.~Semenov.
\newblock Motivic decomposition of a compactification of a {M}erkurjev-{S}uslin
  variety.
\newblock {\em J. Reine Angew. Math.}, 617:153--167, 2008.

\bibitem[Sou85]{Sou-Op-85}
Christophe Soul{\'e}.
\newblock Op\'erations en {$K$}-th\'eorie alg\'ebrique.
\newblock {\em Canad. J. Math.}, 37(3):488--550, 1985.

\bibitem[Sri96]{Sri-96}
Vaseduvan Srinivas.
\newblock {\em Algebraic {$K$}-theory}, volume~90 of {\em Progress in
  Mathematics}.
\newblock Birkh\"auser Boston Inc., Boston, MA, second edition, 1996.

\bibitem[Voe03]{Vo-03}
Vladimir Voevodsky.
\newblock Reduced power operations in motivic cohomology.
\newblock {\em Publ. Math. Inst. Hautes \'Etudes Sci.}, (98):1--57, 2003.

\bibitem[{Zai}10]{Zai-09}
Kirill {Zainoulline}.
\newblock {Degree formula for connective K-theory}.
\newblock {\em Invent. Math.}, 179(3):507--522, 2010.

\end{thebibliography}

\end{document}